\documentclass{amsart}
\usepackage{amscd, amsmath, amsthm, amssymb}

\newtheorem{theorem}{Theorem}[section]
\newtheorem{lemma}[theorem]{Lemma}
\newtheorem{proposition}[theorem]{Proposition}

\theoremstyle{definition}     

\newtheorem{example}[theorem]{Example}

\newtheorem{claim}[theorem]{Claim}

\theoremstyle{remark}

\numberwithin{equation}{section}

\begin{document}

\title[Mordell-Weil group]
{Shioda-Tate formula for an abelian fibered variety and
applications}

\author[K. Oguiso]{Keiji Oguiso}
\thanks {The author was supported by JSPS}
\address{Department of Economics, Keio University, Hiyoshi Kohoku-ku, Yokohama,
Japan and School of Mathematics, Korea Institute for Advanced Study,
Dongdaemun-gu, Seoul 130-722, Korea }
\email{oguiso@hc.cc.keio.ac.jp}

\subjclass[2000]{14C22, 14D06, 14J28, 14J40}

\begin{abstract} We give an explicit formula for the Mordell-Weil rank
of an abelian fibered variety and some of its applications for an
abelian fibered hyperk\"ahler manifold. As a byproduct, we also
give an explicit example of an abelian fibered variety in which
the Picard number of the generic fiber in the sense of scheme is
different from the Picard number of generic closed fibers.
\end{abstract}

\maketitle


\setcounter{section}{0}
\section{Introduction}
The aim of this note is to give an explicit formula for the
Mordell-Weil rank of an abelian fibered variety (Theorem 1.1) and
some of its applications for an abelian fibered hyperk\"ahler
manifold (Theorem 3.1). Our formula is a formal generalization of
the famous formula for an elliptic surface due to Shioda and Tate
([Sh1,2,3]) and some formula for a Calabi-Yau fiber space by
Kawamata [Ka2]. {\it Our formula is also very close to a result of
Bruno Kahn} [Kh], {\it which was recently posted on the ArXiv.}
\par \vskip 1pc

For the precise statement, we introduce a few notations. By an
{\it abelian fibered variety}, we mean a proper surjective
morphism $\varphi : X \longrightarrow Y$ with rational section
$O$, from a normal projective variety $X$ to a normal projective
variety $Y$ such that the generic fiber $A := X_{\eta}$ is a
positive-dimensional abelian variety defined over $\mathbf C(Y)$
(with origin $O$). Here, by a {\it rational section}, we mean a
subvariety $S \subset X$ such that $\varphi \vert S : S
\longrightarrow Y$ is a birational morphism. We assume that:

\begin{list}{}{
\setlength{\leftmargin}{10pt} \setlength{\labelwidth}{6pt} }

\item[(1)] $X$ has only $\mathbf Q$-factorial rational
singularities;

\item[(2)] $Y$ has only $\mathbf Q$-factorial rational
singularities;

\item[(3)] $\varphi$ is equi-dimensional in codimension $1$ in the
sense that there is no prime divisor $D$ on $X$ such that $\dim
\varphi(D) \leq \dim Y - 2$, and

\item[(4)] $h^{1}(X, \mathcal O_{X}) = h^{1}(Y, \mathcal O_{Y})$.
\end{list}

Let $K = \mathbf C(Y)$. The {\it Mordell-Weil group} ${\rm
MW}(\varphi)$ is, by definition, the group $A(K)$ of $K$-rational
points of $A$. In more geometric terms, ${\rm MW}(\varphi)$ is a
group of rational sections of $\varphi$, which we can naturally
regarded as an abelian subgroup of ${\rm Bir}\, (X)$, the group of
the birational transformations of $X$. The N\'eron-Severi group
${\rm NS}\,(A_{K})$ of $A$ is the group of algebraic equivalence
classes of divisors defined over $K$.

Let $\Delta \subset Y$ be the discriminant divisor of $\varphi$
and let $\Delta = \cup_{i=1}^{k} \Delta_{i}$ and
$\varphi^{*}(\Delta_{i})_{\rm red} = \cup_{j=0}^{m_{i}-1} D_{ij}$
be the irreducible decomposition of $\Delta$ and
$\varphi^{*}(\Delta_{i})_{\rm red}$. Note that $m_{i}$ is the
number of prime divisors on $X$ lying over $\Delta_{i}$.
\par \vskip 1pc
Our main result is as follows:

\begin{theorem} \label{theorem:MW}
Under the assumption (1)-(4) and the notation above, the
Mordell-Weil group ${\rm MW}(\varphi)$ is a finitely generated
abelian group of rank
$${\rm rank}\, {\rm MW}\,(\varphi)
= \rho(X) - \rho (Y) - {\rm rank}\, {\rm NS}\,(A_{K}) -
\sum_{i=1}^{k} (m_{i} -1)\,\, .$$ Here $\rho (X)$ (resp. $\rho
(Y)$) is the Picard number of $X$ (resp. of $Y$). In particular,
$${\rm rank}\, {\rm MW}\,(\varphi) \le \rho(X) -2\,\, .$$
\end{theorem}

We believe that our assumptions (1)-(3) are not too restrictive in
the view of minimal model theory in higher dimension. In our
formulation, some condition like (4) is necessary for finite
generation. For instance, in the product case ${\rm pr}_{2} : A
\times Y \longrightarrow Y$, the Mordell-Weil group ${\rm
MW}\,({\rm pr}_{2})$ is far from being finitely generated. {\it
See} [Kh] {\it for a more unified treatment including such cases.}
\par \vskip 1pc

Theorem (1.1) is proved in Section 2. An application for an
abelian fibered hyperk\"ahler manifold will be given in Section 3
(Theorem (3.1)). This part was inspired by a series of work of
Matsushita [Ma1, 2, 3], a work of Sawon [Sa] and also by the author's
previous work [Og]. As one of byproducts of our Theorem
(1.1), we shall give an explicit example of an abelian fibered
hyperk\"ahler manifold $f : X \longrightarrow \mathbf P^{d}$ of
dimension $2d \ge 4$ such that all the smooth fibers $X_{t}$ are isomorphic
to $E_{\zeta}^{d}$, the self-product of the elliptic curve
$E_{\zeta}$ of period $\zeta = e^{2\pi i/3}$, so that $\rho(X_{t}) \ge 2$,
but ${\rm rank}\, {\rm NS} (A_{K}) = 1$ for the generic fiber $A_{K}$
in the sense of scheme (Theorem (3.1)(2)). This example shows that
the Picard number of the generic fiber in the sense of scheme is
not always the same as the Picard number of generic
closed fibers. It might be interesting to compare this with Example
(3.5), a more "moderate" example.
\par \vskip 1pc

{\it Acknowledgement.} An initial idea of this note has been grown
up during my stay at KIAS August 2005 and finalized in the present form
at KIAS March 2007. I would like to express my thanks to Professors
B. Kahn, F. Campana, J. Keum, B. Kim, J.M. Hwang, J. Sawon, D.-Q. Zhang
for several valuable discussions.

\section{Proof of Theorem (1.1)}

In this section, we shall prove Theorem (1.1).
\par \vskip 1pc
By $Z(X)$ (resp. $N(X)$), we denote the free abelian group of Weil
divisors on $X$ (resp. the ablelian group of the linear
equivalence classes of Weil divisors on $X$). By ${\rm Div}\, X$,
we denote the free abelian group of Cartier divisors on $X$. ${\rm
Pic}\, X$ is nothing but the group of linear equivalence classes
of Cartier divisors. Since $X$ is normal and projective, we have a
natural inclusion, ${\rm Pic}\, X \subset N(X)$.
\par \vskip 1pc
Let us start with the following lemma due to Kawamata [Ka1]:
\begin{lemma} \label{lemma:fin}
The quotient group $N(X)/{\rm Pic}\, X$ is a finite abelian group.
\end{lemma}

\begin{proof} Let $\pi
: X' \longrightarrow X$ be a resolution of singularities of $X$.
By [Ka1, Proof
of Lemma (1.1)], ${\rm Pic}\, X'/\pi^{*}{\rm Pic}\, X$ is a
finitely generated
abelian group. Note that ${\rm Pic}\, X' = N(X')$, because $X'$ being smooth.
Thus, the natural surjective homomorphism $N(X') \longrightarrow N(X)$,
together with the projection formula, induces
a surjective homomorphism
${\rm Pic}\, X' / \pi^{*}{\rm Pic}\, X \longrightarrow
N(X)/{\rm Pic}\, X$. Therefore $N(X)/{\rm Pic}\, X$ is also a finitely
generated abelian group. On the other hand, $N(X)/{\rm Pic}\, X$
is a torsion group, because $X$ is $\mathbf Q$-factorial. Combining these two,
we obtain the result.
\end{proof}

\begin{lemma} \label{lemma:eq}
The quotient group ${\rm Pic}\,X/\varphi^{*}{\rm Pic}\, Y$ is a
finitely generated abelian groups of rank $\rho(X) - \rho(Y)$.
\end{lemma}

\begin{proof}
Since $\varphi_{*}\mathcal O_{X} = \mathcal O_{Y}$, the
homomorphism $\varphi^{*}$ is injective. So is its restriction
$(\varphi^{*})^{0} : {\rm Pic}^{0}Y \longrightarrow {\rm Pic}^{0}
X$. Here ${\rm Pic}^{0}Y$ (resp. ${\rm Pic}^{0} X$) is an abelian
variety of dimension $h^{1}(\mathcal O_{Y})$ (resp.
$h^{1}(\mathcal O_{X})$). This is because $Y$ (resp. $X$) has only
rational singularities. Since $h^{1}(\mathcal O_{Y}) =
h^{1}(\mathcal O_{X})$, the homomorphism $(\varphi^{*})^{0}$ is
then an isomorphism. Thus, $\varphi^{*} : {\rm NS}(Y)
\longrightarrow {\rm NS}(X)$ is an injective homomorphism and
${\rm Pic}X/\varphi^{*}{\rm Pic}Y \simeq {\rm
NS}(X)/\varphi^{*}{\rm NS}(Y)$. This implies the result.
\end{proof}

From now, we need some elementary properties of abelian varieties
defined over non-closed fields, for which we refer to the reader
an excellent account of Milne in [AG, Chapter V].
\par \vskip 1pc
Let ${\rm Pic}\, A_{K}$ be the Picard group of $A = A_{K}$, i.e. the
group of linear equivalence classes of divisors defined over $K$,
and let $c : {\rm Pic}\, A_{K} \longrightarrow {\rm NS}\, (A_{K})$
be the natural map. By definition, $c$ is a surjective
homomorphism and its kernel is the group ${\rm Pic}^{0}\, A_{K}$,
i.e. the group of divisors defined over $K$, which are
algebraically equivalent to $0$ modulo linearly equivalent to $0$.
We have ${\rm Pic}^{0}\, A_{K} = \hat{A}_{K}(K)$, where
$\hat{A}_{K}$ is the dual abelian variety of $A = A_{K}$,
both of which are defined over $K$.
\par \vskip 1pc
Let $D$ be a Weil divisor on $X$. Then $D \vert X \setminus {\rm
Sing}\, X$ is Cartier, and therefore, we have a natural surjective
homomorphism:
$$r_{N(X), A} : N(X)\, \longrightarrow {\rm Pic}\, A_{K}\,\, ;\,\,
{\rm cl}(D) \mapsto {\rm cl}(D\vert A)\,\, .$$

Here ${\rm cl}\,(D)$ is the linear equivalence class of $D$. Since
$r_{N(X), A}(\varphi^{*}{\rm Pic}\, Y) = \{0\}$, the homomorphism
$r_{N(X), A}$ induces a homomorphism
$$r_{A} : {\rm Pic}\,X/\varphi^{*}{\rm Pic}\, Y\, \longrightarrow {\rm Pic}\, A_{K}\,\, .$$
By Lemma (2.1), $r_{N(X), A}({\rm Pic}\, X)$ is of finite index in
${\rm Pic}\, A_{K}$. Thus  the image ${\rm Im}\, r_{A}$ is a
subgroup of finite index, of ${\rm Pic}\, A_{K}$. Therefore, by
Lemma (2.2), ${\rm Pic}\, A_{K}$ is also a finitely generated
abelian group of rank $\rho(X) - \rho(Y)$. Now subgroups of ${\rm
Pic}\, A_{K}$ are all finitely generated abelian groups. So, we
may speak of their rank. We denote the $\mathbf Q$-linear
extension of $r_{A}$ by
$$(r_{A})_{\mathbf Q} : ({\rm Pic}\, X/{\rm Pic}\, Y)_{\mathbf Q}
\longrightarrow ({\rm Pic}\, A_{K})_{\mathbf Q}\,\, .$$
$(r_{A})_{\mathbf Q}$ is a surjective linear map.

\begin{lemma} \label{lemma:eqaul}
${\rm MW}(\varphi)$ is a finitely generated abelian group and
satisfies
$${\rm rank}\, {\rm MW}(\varphi) = {\rm rank}\, {\rm Pic}\,
A_{K} - {\rm rank}\, {\rm NS}\, (A_{K})\,\, .$$
\end{lemma}

\begin{proof}
Since $c$ is surjective and its kernel is ${\rm Pic}^{0}\, A_{K}$,
we have:
$${\rm rank}\, {\rm Pic}^{0}\, A_{K} = {\rm rank}\, {\rm Pic}\,
A_{K} - {\rm rank}\, {\rm NS}\, (A_{K})\,\, .$$

As well-known, $\hat{A}_{K}$ and $A_{K}$ are mutually
isogenous over $K$.
Thus, ${\rm MW}(\varphi) = A(K)$ is also a finitely
generated abelian group of the same rank as ${\rm Pic}^{0}\,
A_{K} = \hat{A}_{K}(K)$. Combining this with the formula above,
we obtain the
desired equality.
\end{proof}
\begin{lemma} \label{lemma:neg}

\begin{list}{}{
\setlength{\leftmargin}{10pt} \setlength{\labelwidth}{6pt} }

\item[(1)] $\sum_{i=1}^{k} (m_{i} -1)$
divisor classes
$$D_{ij}\,\, (\,\,1 \le i \le k\,\, , \,\,1 \le j \le m_{i} -1\,\, )$$
are $\mathbf Q$-linearly independent in $({\rm Pic}\, X/\varphi^{*}{\rm Pic Y})_{\mathbf Q}$.

\item[(2)] $\dim\, {\rm Ker}\, ((r_{A})_{\mathbf Q}) = \sum_{i=1}^{k} (m_{i} -1)$.
\end{list}
\end{lemma}

\begin{proof} Let us show (1). There is a large sufficiently divisible number
$M$ such that $MD_{ij}$ are all Cartier. We suffice to show that
there is no non-trivial relation among $1 + \sum_{i=1}^{k}
(m_{i}-1)$ Cartier divisor classes
$$\varphi^{*}L\,\, (L\, \in\, {\rm Pic}\,Y\,),\,\, MD_{ij}\,\,
(\,\, 1 \le i \le k\,\, , \,\,1 \le j \le m_{i} -1\,\, )\,\, .$$
Since $Y$ is projective and $\mathbf Q$-factorial, we may assume
that $L$ is the class of $H_{1} - H_{2}$, where $H_{i}$ are
effective Cartier divisors. Consider the equation in ${\rm Pic}\,
X$ with integral coefficients:

$$c\varphi^{*}H_{1} + \sum_{i = 1}^{k} \sum_{j=1}^{m_{i}-1} c_{ij}D_{ij}
= c\varphi^{*}H_{2}\,\,.$$ Let $C \subset Y$ be a smooth curve
which is a complete intersection of $\dim\, Y -1$ sufficiently
general very ample divisors. Pulling back the above equation to
the induced fibration $X_{C} = X \times_{C} Y \longrightarrow C$,
one obtains the same equation as above, but with $\dim\, C = 1$.
Let $S \subset X_{C}$ be a normal projective surface which is a
complete intersection of $\dim\, X_{C} -2$ sufficiently general
very ample divisors. By pulling back the above equation further to
the induced fibration $S \longrightarrow C$, one obtains the same
equation as above, but now $\dim S = 2$ and $\dim C = 1$. Hence
$c_{ij} = 0$ by the Zariski lemma (see eg. [BHPV, Chap. III, Lemma
(8.3)]). This implies $cL = 0$ in ${\rm Pic}\,(Y)$ as well.
\par \vskip 1pc

Let us show (2). It is clear that the $\mathbf Q$-divisor classes
generated by the divisors in (1) are in the kernel. Let us show
the other inclusion. The argument below is quite close to the
original argument in Shioda-Tate formula (See [Sh1, 3]).
\par \vskip 1pc

Let $D$ be an integral Cartier divisor on $X$ such that ${\rm
cl}(D \vert A) = 0$. It suffices to show that the linear
equivalence class of $D$ is proportional to the sum of divisors in
(1) modulo $\varphi^{*} {\rm Pic}\, Y$. (Indeed, $X$ is $\mathbf
Q$-factorial.)
\par \vskip 1pc

By ${\rm cl}(D \vert A) = 0$, there is an element $f \in K(A)
\setminus \{0\}$ such that $D\vert A = {\rm div} f$ as divisors on
$A$. Since $K(A) = \mathbf C(X)$ and since $K$ is the function
field of $Y$, regarding $f$ as an element of $\mathbf C(X)$, one
can also write this equality as $(D - {\rm div}\, f)\vert A = 0$.
This is again an equality of Weil divisors, but here we regard
${\rm div}\, f$ as a divisor on $X$. Then the support of $D - {\rm
div}\, f$ does not dominate $Y$. So, by replacing $D$ by a
suitable linearly equivalent divisor, we may (and will) assume
that ${\rm Supp}\, D$ does not dominate $Y$. By ${\rm Supp}\, D =
\cup D_{l}$, we denote the irreducible decomposition of ${\rm Supp}\,
D$. Since $\varphi$ is equi-dimensional in codimension $1$ and
$\varphi(D_{l}) \not= Y$, it follows that $\varphi(D_{l})$ are
Weil divisors on $Y$. Again, since both $X$ and $Y$ are $\mathbf
Q$-factorial and since $\varphi$ is equi-dimensional in
codimension $1$, there are positive integers $n_{l}$ and $n_{ij}'$
such that each $n_{l}D_{l}$ is either an element of $\varphi^{*}
{\rm Pic}\,Y$, or one of $n_{ij}'D_{ij}$ ($1 \le i \le k$, $0 \le
j \le m_{i} - 1$). Since $\sum_{j=0}^{m_{i}-1} c_{ij}D_{ij} \in
\varphi^{*}{\rm Pic}\, Y$ for suitably chosen $c_{ij} > 0$, the
result follows.
\end{proof}

By Lemma (2.4), we have
$${\rm rank}\,{\rm Pic}\, A_{K} = \rho(X) - \rho(Y) - \sum_{i=1}^{k}
(m_{i} -1)\,\, .$$ Substituting this into the equality in Lemma
(2.3), we obtain the desired formula in Theorem (1.1).

\section{Applications for an abelian fibered hyperk\"ahler manifold}

In this section, we shall work in the category of complex analytic
spaces. Especially, unless stated otherwise,
the term {\it generic point} means a closed
point being in the complement of the union of countably many
proper closed analytic subvarieties.
\par \vskip 1pc
By a {\it hyperk\"ahler manifold} we mean a compact
simply-connected complex K\"ahler manifold $M$ having everywhere
non-degenerate holomorphic $2$-form $\sigma_{M}$ such that
$H^{0}(M, \Omega_{M}^{2}) = \mathbf C \sigma_{M}$. We note that
$\dim_{\mathbf C}\, M$ is even.
\par \vskip 1pc
For example, a K3 surface $S$ is a $2$-dimensional hyperk\"ahler
manifold and its Hilbert scheme of $0$-dimensional closed
subscheme of length $d$
$$S^{[d]} := {\rm Hilb}^{d}\, S$$
is a hyperk\"ahler manifold of dimension $2d$ ([Be]).
\par \vskip 1pc
By a result of Matsushita [Ma1,2], any abelian fibered
hyperk\"ahler manifold satisfies the conditions (1)-(4) required
in Theorem (1.1) except perhaps the projectivity. See also
Proposition (3.2) about the projectivity of a fibered
hyperk\"ahler manifold.
\par \vskip 1pc
We refer to the readers an excellent account [GHJ, Part III] for
basic properties about hyperk\"ahler manifolds.
\par \vskip 1pc
In this section, we shall show the following:

\begin{theorem} \label{theorem:hk}

(1) For each positive integer $d \ge 2$ and each integer $\rho$
with $2 \le \rho \le 20$, there is an abelian fibered
hyperk\"ahler manifold $f_{d, \rho} : M_{d, \rho} \longrightarrow
\mathbf P^{d}$ of dimension $2d$ such that $\rho(M_{d, \rho}) =
\rho$ and ${\rm rank}\, {\rm MW}(f_{d, \rho}) = \rho - 2$.

(2) For each positive integer $d \ge 2$, there is an abelian
fibered hyperk\"ahler manifold $f_{d} : M_{d} \longrightarrow
\mathbf P^{d}$ of dimension $2d$ such that the Picard number of
all smooth closed fibers are greater than $1$ but the Picard
number of the generic fiber over $\mathbf C(\mathbf P^{d})$ is
$1$. More geometrically, all the smooth fibers are isomorphic to
$E_{\zeta}^{d}$ but the generic fiber
$A = M_{\eta}$ is simple over $\mathbf C(\mathbf P^{d})$. Here
$E_{\zeta}$ is the elliptic curve of period $\zeta = e^{2\pi
i/3}$ and $E_{\zeta}^{d}$ is the $d$-th self-product of $E_{\zeta}$.
\end{theorem}

In the rest, we shall prove Theorem (3.1).
\par \vskip 1pc
Let us start with the following projectivity criterion due to F.
Campana:

\begin{proposition} \label{proposition:proj}
Let $f : M \longrightarrow B$ be a surjective morphism from a
hyperk\"ahler manifold to a normal projective variety $B$ with
connected fibers. Assume that $0 < \dim\, B < \dim\, M = 2d$ and $M$
has a subvariety $S$ such that $f(S) = B$ and $\dim\, S = \dim\,
B$. Then $M$ is projective and the generic fiber $M_{t}$ is an
abelian variety.
\end{proposition}

\begin{proof} By definition, $M$ is K\"ahler. By Matsushita [Ma1],
any generic fiber $M_{t} = f^{-1}(t)$ is a smooth Lagrangian
submanifold, i.e. $\sigma_{M} \vert M_{t} = 0$ and $\dim\, M_{t} =
\dim\, M/2 = d$. (Note that in [Ma1],
$M$ is also assumed to be projective. However, the argument there
works for non-projective $M$ if one uses a K\"ahler class of $M$
instead of an ample class of $M$ used there.) This implies that
$$\Omega^{1}_{M_{t}} \simeq N_{M_{t}/M} \simeq
\mathcal O_{M_{t}}^{\oplus d}\,$$ and therefore, that $M_{t}$ is a
$d$-dimensional complex torus. Moreover, by Voisin's lemma [Ca3],
$M_{t}$ is also projective. Thus $M_{t}$ is an abelian variety.
Since $B$ is projective and $f(S) = B$, we have
$$d = \dim\, B = a(B) \le a(S) \le \dim\, S = d\,\, .$$
Here $a(S)$ (resp. $a(B)$) is the algebraic dimension of $S$
(resp. $B$). Thus, $a(S) = \dim\, S = d$. Hence $S$ is
bimeromorphic to a projective variety. In particular, $S$ is
covered by complete algebraic curves. Thus, any two generic points
in $M$ can be connected by a chain of complete algebraic curves.
Therefore, $M$ is Moishezon, i.e. $\dim\, M = a(M)$, by a result
of Campana [Ca1] (See also [Ca2] for a different approach). Since
$M$ is K\"ahler, $M$ is then projective by Moishezon's criterion
[Mo].
\end{proof}

\begin{proposition} \label{proposition:deform}
Assume that there is an abelian fibered (projective) hyperk\"ahler
manifold $f : M \longrightarrow \mathbf P^{d}$ having holomorphic
zero section $O$ and $r$ holomorphic sections, say $S_{i}$ ($1 \le
i \le r)$, which are linearly independent in ${\rm MW}(f)$. Then,
there is an abelian fibered projective hyperk\"ahler manifold $f'
: M' \longrightarrow \mathbf P^{d}$ such that $\dim\, M' = \dim\,
M$, ${\rm rank}\, {\rm MW}(f') = r$, and $\rho(M') = r + 2$.
\end{proposition}

In general, by a {\it holomorphic section} of a fiber space
$\varphi : X \longrightarrow Y$, we mean a subvariety $S \subset
X$ such that $\varphi \vert S : S \longrightarrow Y$ is an
isomorphism. One can then identify $S$ with a holomorphic map $Y
\longrightarrow X$ given by $t \mapsto (\varphi \vert S)^{-1}(t)$.

\begin{proof} Let $u : \mathcal U \longrightarrow \mathcal K$
be the Kuranishi family of $M$, in which we assume that $0 \in
\mathcal K$, $M = \mathcal U_{0} := u^{-1}(0)$. Each fiber
$\mathcal U_{t} := u^{-1}(t)$ ($t \in \mathcal K$) is a
hyperk\"ahler manifold. Let $\sigma_{t}$ be an everywhere
non-degenerate holomorphic $2$-form on $\mathcal U_{t}$. Then, by
choosing a marking $\iota : R^{2}u_{*} \mathbf Z \simeq \Lambda
\times \mathcal K$, we can define the period map $p : \mathcal K
\longrightarrow \mathcal P\, ;\, t \mapsto [\iota(\sigma_{t})]$.
Here
$$\mathcal P := \{[\sigma] \in \mathbf P(\Lambda_{\mathbf
C})\, \vert\, (\sigma, \sigma) = 0\,\, ,\,\,  (\sigma,
\overline{\sigma})
> 0\,\}\,\, \subset \mathbf P(\Lambda_{\mathbf C})\,\, .$$
By the local Torelli theorem, $p$ is a local isomorphism (See eg.
[JHG, Part III]). Thus we can (and will) identify $\mathcal K$
with a small open neighborhood of $p(0)$ in $\mathcal P$, say $W$,
via $p$.
\par \vskip 1pc
{\it In what follows, we freely shrink $0 = p(0) \in W$ whenever we need
to do so.}
\par \vskip 1pc
Let $L$ be the pullback of the hyperplane on $\mathbf P^{d}$ by
$f$. Then $f = \Phi_{\vert D \vert}$. Here $\Phi_{\vert D \vert}$
is the morphism associated with the complete linear system $\vert
L \vert$. Note that $O$ and $S_{i}$ are Lagrangian submanifolds of
$M$, because they are isomorphic to $\mathbf P^{d}$ and $\mathbf
P^{d}$ has no non-zero global holomorphic $2$-form. Note also that
$H^{2}(\mathbf P^{d}, \mathbf Z) \simeq \mathbf Z$.

Let us consider the deformation of $M$, say $W(r) (\subset W)$,
which keeps the class $[L]$ being $(1,1)$-class and the Lagrangian
submanifolds $O$, $S_{i}$ ($1 \le i \le r$) being Lagrangian. As
in Sawon [Sa], it follows from [Ma2] and [Vo] that $W(r)$ is the
intersection of $W (\subset \mathbf P(\Lambda_{\mathbf C}))$ with
$2 + r$ rational hyperplanes, say $H_{j}$ ($1 \le j \le 2+r$),
corresponding to the required conditions on $L$, $O$ and $S_{i}$
($1 \le i \le r$). Note that $0 \in W(r)$.

Let $w : \mathcal M \longrightarrow W(r)$ be the family induced
from the universal family $u : \mathcal U \longrightarrow \mathcal
K = W$. By construction, the line bundle $L$ and Lagrangian
submanifolds $O$ and $S_{i}$ ($1 \le i \le r$) are all extended
over $W(r)$, say, $\mathcal L$, $\mathcal O$ and $\mathcal S_{i}$.
It is shown by [Ma3] (see also [Sa]) that $h^{0}(\mathcal M_{t},
\mathcal L_{t}) = d+1$ and $h^{i}(\mathcal M_{t}, \mathcal L_{t})
= 0$ ($i \ge 1$) for all $t \in W(r)$. Thus, $w_{*}\mathcal L$ is
a locally free sheaf and satisfies the base change property. Thus,
by the freeness of $\vert L \vert$, the linear system $\vert
\mathcal L \vert$ is $w$-free and gives a fibration $\tilde{f} :
\mathcal M \longrightarrow \mathbf P^{d} \times W(r)$ with
holomorphic sections $\mathcal O$, $\mathcal S_{i}$ ($1 \le i \le
r$), over $W(r)$. Here, we also used the fact that $\mathbf P^{d}$
is rigid.

Fiberwisely, $\tilde{f}$ induces an abelian fibration
$\tilde{f}_{t} : \mathcal M_{t} \longrightarrow \mathbf P^{d}$ by
[Ma1](see also Proposition (3.2)) with holomorphic sections
$\mathcal O_{t} := \mathcal O \vert \mathcal M_{t}$ and $\mathcal
S_{i, t} := \mathcal S_{i} \vert \mathcal M_{t}$ ($1 \le i \le
r$), for each $t \in W(r)$. By Proposition (3.2), $\mathcal M_{t}$
($t \in W(r)$) is projective. We take $\mathcal O_{t}$ as the
origin of the Mordell-Weil group ${\rm MW}(\tilde{f}_{t})$. The
sections $\mathcal S_{i, t}$ ($1 \le i \le r$) are naturally
regarded as elements of ${\rm MW}(\tilde{f}_{t})$.

\begin{claim}\label{claim:down} For a generic point $t \in W(r)$, one has:
\begin{list}{}{ \setlength{\leftmargin}{10pt}
\setlength{\labelwidth}{6pt} }

\item[(1)] $\rho(\mathcal M_{t}) \le 2 + r$ and

\item[(2)] $\mathcal S_{i, t}$ ($1 \le i \le r$) are linearly
independent in the Mordell-Weil group ${\rm MW}\,
(\tilde{f}_{t})$.
\end{list}
\end{claim}
\begin{proof} By definition of generic point, we may (and will) show (1)
and (2) individually.
\par \vskip 1pc
Let us show (1). Let $l$ be the number of independent rational
hyperplanes among $H_{j}$ ($1 \le j \le 2 + r$). Let $V \subset
W(r)$ be the set of points $t \in W(r)$ such that the number of
independent element $h \in \Lambda$ satisfying that $(h,
p(\sigma_{t})) = 0$ is exactly $l$.

Since each rational hyperplane in $\mathbf P(\Lambda_{\mathbf C})$
is of the form $(h, *) = 0$ for some $h \in \Lambda \setminus
\{0\}$, they are countable in number.  Thus $V$ is the complement
of the union of countably many proper closed analytic subsets of
$W(r)$.

Let $t \in V$. Then, by the Lefschetz $(1,1)$-Theorem,
$\rho(\mathcal M_{t}) = l \le 2 + r$. This implies (1).
\par \vskip 1pc

Let us show (2). Let $c := (c_{i})_{i=1}^{r} \in \mathbf Z^{r}
\setminus \{0\}$. Consider the subset
$$\mathcal D_{c} := \{\, t \in W(r)\, \vert\,
\sum_{i=1}^{r} c_{i}\mathcal S_{i, t} = \mathcal O_{t}\,\, {\rm
in}\,\, {\rm MW}(\tilde{f}_{t})\, \}\,\, ,$$
and its Zariski closure $\overline{\mathcal D}_{c}$ in $W(r)$.
Then $\overline{\mathcal D}_{c}$ is a closed analytic subset of
$W(r)$.

Since $\mathbf Z^{r} \setminus \{0\}$ is a countable set, it now
suffices to show that $0 \not\in \overline{\mathcal D}_{c}$ for
each $c$.

Choose a fiber $F$ of $f = \tilde{f}_{0} : M \longrightarrow
\mathbf P^{d}$ such that $F$ is an abelian variety, with origin
$\mathcal O \vert F$, in which $\mathcal S_{i}\vert F (= S_{i}
\vert F)$ ($1 \le i \le r$) are linearly independent points. By
assumption, such a fiber $F$ exists. Choose then a smooth
subfamily of $\tilde{f}$, say,
$$\pi := \tilde{f}_{\mathcal F} :
\mathcal F \longrightarrow W(r)$$ such that $F = \mathcal F_{0} =
\pi^{-1}(0)$ and $\mathcal F_{t} = \pi^{-1}(t)$ is a fiber of
$\tilde{f}_{t}$, which is an abelian variety with origin $\mathcal
O \vert \mathcal F_{t}$. Again such a family exists. Note also
that $\mathcal O \vert \mathcal F$ and $\mathcal S_{i} \vert
\mathcal F$ are holomorphic sections of $\pi$. Since an abelian
variety contains no rational curve, the sum $\sum_{i=1}^{r}
c_{i}\mathcal S_{i} \vert \mathcal F$, as well as $\mathcal O
\vert \mathcal F$ and $\mathcal S_{i} \vert \mathcal F$, defines
(not only a rational section but also) a holomorphic section of
$\pi$.

Regard $\mathcal O \vert \mathcal F$ and $\sum_{i=1}^{r}
c_{i}\mathcal S_{i} \vert \mathcal F$ as holomorphic maps from
$W(r)$ to $\mathcal F$, rather than subvarieties of $\mathcal F$.
Then, by definition of $\mathcal D_{c}$, one has:
$$(\sum_{i=1}^{r} c_{i}\mathcal S_{i} \vert \mathcal F)(t)
= (\mathcal O \vert \mathcal F)(t) \,\, $$ for all $t \in \mathcal
D_{c}$. Since the sections $\sum_{i=1}^{r} c_{i}\mathcal S_{i}
\vert \mathcal F$ and $\mathcal O \vert \mathcal F$ are both
holomorphic, we obtain
$$(\sum_{i=1}^{r} c_{i}\mathcal S_{i} \vert \mathcal F)(t')
= (\mathcal O \vert \mathcal F)(t')\,\,$$ for all $t' \in
\overline{\mathcal D}_{c}$. Thus $0 \not\in \overline{\mathcal
D}_{c}$. Indeed, otherwise, we would have
$$(\sum_{i=1}^{r} c_{i}\mathcal S_{i} \vert \mathcal F)(0)
= (\mathcal O \vert \mathcal F)(0)\,\, ,\,\,{\rm i.e.}\,\,
\sum_{i=1}^{r} c_{i}S_{i} \vert F = O \vert F\,\, ,$$ a
contradiction to the choice of $F$.
\end{proof}

Let us return back to the proof of Proposition (3.3). Let $t \in
W(r)$ be a generic point in Claim (3.4). Then ${\rm rank}\, {\rm
MW}(f_{t}) \ge r$. Hence $\rho(\mathcal M_{t}) \ge 2 + r$ by
Theorem (1.1). Combining this with Claim (3.4)(1), we obtain
$\rho(\mathcal M_{t}) = 2 + r$, and therefore, ${\rm MW}(f_{t}) =
r$. We may now take this $\tilde{f}_{t} : \mathcal M_{t}
\longrightarrow \mathbf P^{d}$ as $f' : M \longrightarrow \mathbf
P^{d}$.
\end{proof}

Now we are ready to prove Theorem (3.1).
\par \vskip 1pc
{\it Proof of Theorem (3.1)(1).}
\par \vskip 1pc
 As well-known (see eg. [Og]),
there is an elliptic K3 surface $\varphi : S \longrightarrow
\mathbf P^{1}$ with section $O$ such that $\rho(S) = 20$ and ${\rm
rank}\, {\rm MW}\, (\varphi) = 18$.

We have an abelian fibration
$$\varphi_{d} : S^{[d]} \longrightarrow
\mathbf P^{d}\, =\, {\rm Sym}^{d}\, \mathbf P^{1}\,\, ,$$ which is
the composition of the Hilbert-Chow morphism $S^{[d]}
\longrightarrow {\rm Sym}^{d}\, S$ and the natural map ${\rm
Sym}^{d}\, S \longrightarrow {\rm Sym}^{d}\, \mathbf P^{1}$
induced by $\varphi$. Moreover, the $0$-section and $18$
independent sections of $\varphi$ give rise to the holomorphic
$0$-section and $18$ independent holomorphic sections of
$\varphi_{d}$. Here we also used the fact that ${\rm Hilb}^{d}\,
\mathbf P^{1} \simeq \mathbf P^{d}$ and there is no non-trivial
birational morphism from $\mathbf P^{d}$ to $\mathbf P^{d}$.

Applying Proposition (3.3) to this $\varphi_{d}$, we obtain
Theorem (3.1)(1).

\par \vskip 1pc
{\it Proof of Theorem (3.1)(2).}
\par \vskip 1pc

Let $f : S \longrightarrow \mathbf P^{1}$ be an elliptic K3
surface defined by the Weierstrass equation
$$y^{2} = x^{3} - (t^{11} -1)\,\, .$$
It is well-known that $\rho(S) = 2$ (see eg. [Ko]). This follows
from the fact that $S$ admits a non-symplectic automorphism $g$ of
maximum order $66$:
$$g^{*}(x, y,t) = (\zeta_{3} x, -y, \zeta_{11} t)\,\, .$$
Here $\zeta_{n} = e^{2\pi i/n}$. Indeed, $g$ then acts on the
space of global holomorphic $2$-forms as
$$g^{*}\frac{dx \wedge dt}{y} =
\zeta_{66}^{-5}\frac{dx \wedge dt}{y}\,\, .$$ This implies that
$$20 = \varphi(66)\, \vert\, {\rm rank}\, T(S) = 22 - \rho(S) \le 21\,\, .$$
Here $\varphi$ is the Euler function and $T(S)$ is the
transcendental lattice of $S$ (cf. [BHPW]). Thus, $\rho(S) = 2$.
From the Weierstrass equation, we also see that each smooth fiber
of $f$ is isomorphic to $E_{\zeta}$ and the singular fibers are
all isomorphic to the cuspidal rational curve. Let $f_{d} :
S^{[d]} \longrightarrow \mathbf P^{d}$ be the abelian fibration
induced by $f$. Each smooth fiber of $f_{d}$ is then isomorphic to
the product abelian variety $E_{\zeta}^{d}$.

Let us consider the generic fiber $A$ of $f_{d}$ in the sense of
scheme. $A$ is an abelian variety defined over $K = \mathbf
C(\mathbf P^{d})$. We shall compute the rank of ${\rm NS}(A_{K})$
by using Theorem (1.1).

By $\rho(S) = 2$, we have
$$\rho(S^{[d]}) = \rho(S) + 1 = 3\,\, .$$
Here "$+1$" comes from the exceptional divisor $E$ of the
Hilbert-Chow morphism, which is generically the blow-up of the big
diagonal of the Chow variety ${\rm Sym}^{d}\, S$. Observe that the
divisor $f_{d}^{-1}(\Delta)$, where $\Delta$ is the big diagonal
of $\mathbf P^{d} = {\rm Sym}^{d}\, \mathbf P^{1}$, consists of
two irreducible components. Indeed, one component is $E$ and the
other component is the proper transform of
$\overline{f}_{d}^{-1}(\Delta)$, where $\overline{f}_{d}$ is the
natural morphism from ${\rm Sym}^{d}\, S$ to ${\rm Sym}^{d}\,
\mathbf P^{1}$. Then, by Theorem (1.1),
$$0 \le {\rm rank}\, {\rm MW}\, (f_{d}) \le 3 - 1 - 1 - {\rm
rank}\, {\rm NS}\, (A_{K}) = 1 - {\rm rank}\, {\rm NS}\,
(A_{K})\,\, .$$ Since ${\rm rank}\, {\rm NS}(A_{K}) \ge 1$, this
implies that ${\rm rank}\, {\rm NS}\, (A_{K}) = 1$. In particular,
$A$ is simple, i.e. $A$ is not isogenous to the product of lower
dimensional abelian varieties, over $\mathbf C(\mathbf P^{d})$.
\par \vskip 1pc
This completes the proof of Theorem (3.1).
\par \vskip 1pc

It may be interesting to compare Theorem (3.1)(2) with the
following probably more standard:
\begin{example}\label{example:schoen} Let $f_{i} : S_{i} \longrightarrow
\mathbf P^{1}$ ($i = 1$, $2$) be two relatively minimal rational
elliptic surfaces with section, such that discriminant locus
$\Delta_{i} \subset
\mathbf P^{1}$ of $f_{i}$ are disjoint. Then, as is shown by [Sc],
$X := S_{1} \times_{\mathbf P^{1}} S_{2}$ is a smooth Calabi-Yau
threefold having an abelian fibration $f : X \longrightarrow
\mathbf P^{1}$ induced by $f_{1}$ and $f_{2}$. In this example,
assume further that $f_{i}$ are generic in the sense that all
singular fibers are nodal rational curves and that the generic
fiber $S_{i, \eta}$ ($i = 1$, $2$) are not mutually isogenous.
(Most cases are such cases.)

Then the Picard number of any closed generic fiber $S_{1, t}
\times S_{2, t}$ and the Picard number of the generic fiber $S_{1,
\eta} \times_{{\rm Spec}\, \mathbf C(\mathbf P^{1})} S_{2, \eta}$
in the sense of scheme are all $2$. This follows from the assumption
that two direct factors are not mutually isogenous. (cf. [ibid]).

Under the same assumption, by using Theorem (1.1), we also obtain
that:
$${\rm rank}\,{\rm MW}\,(f)\, =\, 19 - 1 - 2\, =\, 16\, =\,
\,8 + 8\, =\, {\rm rank}\,{\rm MW}\,(f_1) + {\rm rank}\,{\rm
MW}\,(f_2)\,\, .$$
\end{example}

\end{document}